\documentclass{amsart}%[dvipdfmx]

\makeatletter
	
	\@addtoreset{equation}{section}

	\@addtoreset{figure}{section}
\makeatother

\usepackage{amsmath}
\usepackage{amsfonts}
\usepackage{graphicx}
\usepackage{subcaption}
\usepackage{url}
\usepackage{mathtools}
\usepackage{caption}
\usepackage{appendix}

\usepackage[driver=dvipdfmx,text={372pt,533pt},centering]{geometry}

\numberwithin{equation}{section}

\newtheorem{thm}{Theorem}[section]
\newtheorem{defi}[thm]{Definition}%[section]
\newtheorem{prop}[thm]{Proposition}%[section]
\newtheorem{cor}[thm]{Corollary}
\newtheorem{conj}[thm]{Conjecture}
\newtheorem{ques}[thm]{Question}

\newtheorem{rem}[thm]{Remark}

\newcommand{\iu}{\sqrt{-1}}
\newcommand{\id}{\mathrm{id}}
\newcommand{\A}{\mathcal{A}}

\newcommand{\bb}{\boldsymbol{b}}
\newcommand{\bk}{\boldsymbol{k}}
\newcommand{\bp}{\boldsymbol{p}}
\newcommand{\bw}{\boldsymbol{w}}
\DeclareMathOperator{\tr}{Tr}
\DeclareMathOperator{\Vol}{Vol}
\DeclareMathOperator{\li}{Li_2}
\DeclareMathOperator{\im}{Im}

\DeclareMathOperator{\sgn}{sgn}

\title[]{On the Potential Function of the Colored Jones Polynomial and the AJ conjecture}%The Potential Function and the AJ conjecture
\author[]{Shun Sawabe}
\address{Department of Pure and Applied Mathematics,\\ School of Fundamental Science and Engineering,\\
Waseda University,\\
3-4-1 Okubo, Shinjuku, Tokyo 169-8555, Japan}
\email{ssawabe[at]aoni.waseda.jp}
\date{}

\begin{document}

\maketitle

\begin{abstract}
	The AJ conjecture states that the $A_q$-polynomial, which derives from a recurrence relation of the colored Jones polynomial,
	evaluated at $q=1$ essentially equals the $A$-polynomial.
	We can also obtain a factor of the $A$-polynomial from the potential function of the colored Jones polynomial.
	In this paper, we introduce a polynomial obtained from a potential function of the colored Jones polynomial based on those studies,
	and show that the polynomial appears as a factor of the $A_q$-polynomial evaluated at $q=1$. 
\end{abstract}

\section{Introduction}
Quantum invariants are closely related to $3$-dimensional geometry.
One of the crucial conjectures is the volume conjecture.
The volume conjecture states that a certain limit of the colored Jones polynomial for a knot $K$
evaluated at the root of unity equals the volume of the knot complement.
Note that in the statement of the volume conjecture, we have to normalize the colored Jones polynomial
so that the value for the unknot $ \bigcirc $ is $1$.
To distinguish the normalization, we write $J_n(K;q)$ to denote the $n$-th colored Jones polynomial satisfying
$J_n(\bigcirc;q) = [n]$, where
\[
	[n] = \frac{\{n\}}{\{1\}},\quad \{n\} = q^{\frac{n}{2}}-q^{-\frac{n}{2}}
\]
for an integer $n$, and $J'_n(K;q)$ for the one satisfying $J'_n(\bigcirc;q) = 1$ in this paper.
\begin{conj}[Volume Conjecture \cite{MM}]
	For any knot $K$, the colored Jones polynomial $J'_N(K;q)$ satisfies
	\[
		2 \pi \lim _{N \to \infty} \frac{\log |J'_N(K;q=\xi_N)|}{N} = v_3 ||K||,
	\]
	where $\xi_N = e^{\frac{2 \pi \iu}{N}}$,
	$v_3$ is the volume of the ideal regular tetrahedron in the three-dimensional hyperbolic space
	and $|| \cdot ||$ is the simplicial volume for the complement of the knot $K$.
\end{conj}
Yokota \cite{Yo} gave the idea to prove this conjecture.
In his theory, a potential function plays an important role.
\begin{defi}
	Suppose that the asymptotic behavior of a certain quantity $Q_N$ for a sufficiently large $N$ is %can be written as
	\[
		Q_N \sim \int \cdots \int _\Omega P_N e^{\frac{N}{2 \pi \iu} \Phi(z_1,\ldots,z_\nu)}dz_1 \cdots dz_\nu,
	\]
	where $P_N$ grows at most polynomially and $ \Omega $ is a region in $\mathbb{C}^\nu$. We call this function $\Phi(z_1,\ldots,z_\nu)$ a potential function of $Q_N$.
\end{defi}
Yokota considered the potential function of the Kashaev invariant
and established a relationship between a saddle point equation and a triangulation of a knot complement.
We also considered the potential function of the colored Jones polynomial with parameters
corresponding to the colors and found a geometric meaning of derivatives with those parameters in \cite{Sa}.
In the knot case, the upshot is as follows:
Let $n$ be an integer, let $m$ be a half-integer satisfying $n=2m+1$,
and let $K$ be a hyperbolic knot with a diagram $D$.
We also let $F^D(m,k_1,\ldots,k_\rho)$ be the summand of the colored Jones polynomial
obtained by assigning the $R$-matrix to crossings of the diagram $D$.
Let $\Phi_D(\alpha,w_1,\ldots,w_\rho)$ be a potential function of the colored Jones polynomial
$J_n(K;\xi_N)$ constructed from the summand $F^D(m,k_1,\ldots,k_\rho)$.
In this process, the parameter $\alpha$ corresponds to $m$,
and the parameters $w_1,\ldots,w_\rho$ respectively correspond to $k_1,\ldots,k_\rho$.
Then, %the system of equations
\begin{equation} \label{eq:phihyp}
	\begin{dcases}
		\exp \left(w_i \frac{\partial \Phi _D}{\partial w_i} \right)= 1, &(i = 1,\ldots,\rho)\\
		\exp \left(\alpha \frac{\partial \Phi _D}{\partial \alpha}\right) = l^2
	\end{dcases}
\end{equation}
is a necessary condition that the knot complement admits a hyperbolic structure.
Here, $l^2$ is the dilation component of the preferred longitude.
We use this notation so that it matches the variable of the $A$-polynomial later.
In the above explanation, we use the summand obtained directly by the
$R$-matrix formulation of the colored Jones polynomial.
The colored Jones polynomial, however, has many presentations. For example,
we often use a $q$-analog of the binomial theorem to simplify the formula.
In general, a summand $F(m,k_1,\ldots,k_\nu)$ of the colored Jones polynomial is of the form
\[
	F(m,k_1,\ldots,k_\nu)=(-1)^{(\boldsymbol{e},\tilde{\bk})}q^{\frac{1}{2}\tilde{\bk}^\top A \tilde{\bk}+(\bb,\tilde{\bk})+c}\frac{\prod_i (q)_{(\bp_i,\tilde{\bk})+r_i}}{\prod_j (q)_{(\bp'_j,\tilde{\bk})+r'_j}},
\]
where $\boldsymbol{e}$ is a $(\nu+1)$-vector in $\{0,1\}^{\nu+1} \subset \mathbb{Q}^{\nu+1}$,
$ \tilde{\bk}=(m,k_1,\ldots,k_\nu)$,
$A$ is a $(\nu+1)$-th symmetric matrix with integer entries,
$\bb$ is a vector in $\mathbb{Q}^{\nu+1}$, $c$ is an integer, $\bp_i$ and $\bp'_j$ are $(\nu+1)$-vectors in $\mathbb{Z}^{\nu+1}\subset \mathbb{Q}^{\nu+1}$,
and $r_i$, $r'_j$ are integers.
We can construct a potential function $\Phi_F(\alpha,w_1,\ldots,w_\nu)$ from the summand $F(m,k_1,\ldots,k_\nu)$, and the equations
\begin{equation} \label{eq:delphiFitr}
	\begin{dcases}
		\exp \left(w_i \frac{\partial \Phi_F}{\partial w_i} \right)&= 1, (i = 1,\ldots,\nu)\\
		\exp \left(\alpha \frac{\partial \Phi_F}{\partial \alpha}\right) &= l^2,
	\end{dcases}
\end{equation}
are equivalent to algebraic equations
\[
	p_0(l,w_1,\ldots,w_\nu) = 0,\ p_1(w_1,\ldots,w_\nu) = 0,\ldots,\ p_{\nu}(w_1,\ldots,w_\nu) = 0,
\]
where $p_0(l,w_1,\ldots,w_\nu) \in \mathbb{Q}(\alpha)[l,w_1,\ldots,w_\nu]$, and
$p_1(w_1,\ldots,w_\nu),\ldots,p_\nu(w_1,\ldots,w_\nu)$ are polynomials in $\mathbb{Q}(\alpha)[w_1,\ldots,w_\nu]$.
We suppose that those equations have solutions.
\begin{defi}
	We say that a summand $F(m,k_1,\ldots,k_\nu)$ is solvable if the system of equations \eqref{eq:delphiFitr} has solutions.
\end{defi}
At least, the summand obtained directly from a knot diagram is solvable
because there exists a solution that gives a hyperbolic structure.
We eliminate variables $w_1,\ldots,w_\nu$ from here.
\begin{defi}
	We define the polynomial $a_F(l,\alpha)$ as the generator of the ideal
	\[
		I = \langle p_0,p_1,\ldots,p_\nu \rangle \cap \mathbb{Q}(\alpha)[l] \subset \mathbb{Q}(\alpha)[l]
	\]
	whose coefficients are in $\mathbb{Z}[\alpha]$ and coprime.
	The polynomial $a_F(l,\alpha)$ is defined up to multiplication by $ \pm 1$.
\end{defi}
Considering that the system of equations \eqref{eq:phihyp} is a necessary condition that the knot complement admits a hyperbolic structure \cite{Sa},
and that the factor of the $A$-polynomial corresponding to irreducible representations is obtained from
a potential function by eliminating redundant variables \cite{Hi,Y2}, the following question is reasonable:
\begin{ques}
	Does the polynomial $a_F(l,\alpha)$ coincide with the nonabelian $A$-polynomial (up to multiplication by $ \pm 1$)
	for any solvable summand $F(m,k_1,\ldots,k_\nu)$ of the colored Jones polynomial?
\end{ques}
See also \cite{Gu} for such relation between the colored Jones polynomial and the $A$-polynomial.
On the other hand, the AJ conjecture is known as a relationship between the colored Jones polynomial and the $A$-polynomial.
The AJ conjecture states that the $A$-polynomial is obtained from the recurrence relation of the colored Jones polynomial
by evaluating $q$ at $1$.
In other words, the recurrence relation of the colored Jones polynomial $J_K(n) = J_n(K;q)$ for a knot $K$ determines a noncommutative polynomial $A_q(K)(E,Q)$,
and this polynomial is conjectured to be a $q$-version of the $A$-polynomial.
Here, $E$ and $Q$ are operators defined by
\[
	(EJ_K)(n) = J_K(n+1),\quad (QJ_K)(n) = q^n J_K(n).
\]
\begin{conj}[the AJ conjecture \cite{Ga}]
	For any knot $K$, the $A$-polynomial $A_K(l,\alpha)$ is equal to the $A_q$-polynomial $ \varepsilon A_q(K)(l,\alpha^2)$
	up to multiplication by an element in $ \mathbb{Q}(\alpha)$,
	where $ \varepsilon $ is an evaluation map at $q=1$.
\end{conj}
Garoufalidis \cite{Ga} proposed this conjecture.
The AJ conjecture holds for twist knots \cite{Le}, torus knots \cite{Tr}, and
some classes of hyperbolic knots \cite{LZ}.
Detcherry and Garoufalidis \cite{DG} also considered the AJ conjecture from the perspective
of a triangulation of the knot complement.
Takata \cite{Ta} observed it with twist knots and discussed the results in terms of the volume conjecture. \par
In this paper, we connect these conjectures on the relationship
between the colored Jones polynomial and the $A$-polynomial via the potential function.
The $A_q$-polynomial is obtained from the summand of the colored Jones polynomial by creative telescoping (See \cite{Ga,GL} or Section \ref{sec:AJ}).
We will compare this process with the above conjectural method to obtain the $A$-polynomial.
For an integer $n$, let $F=F(m,k_1,\ldots,k_\nu)$ be the summand of the colored Jones polynomial $J_K(n)=J_n(K;q)$ for a knot $K$,
where $m$ is a half-integer satisfying $n=2m+1$.
Let $E,\ E_{(0)},\ E_j,\ Q,\ Q_{(0)}$, and $Q_j$ be operators defined by
\begin{align*}
	(EF) \left(\frac{n-1}{2},k_1,\ldots,k_\nu \right) &= F \left(\frac{n}{2},k_1,\ldots,k_\nu \right),\\
	(E_{(0)}F)(m,k_1,\ldots,k_\nu) &= F(m+1,k_1,\ldots,k_\nu),\\
	(E_jF)(m,k_1,\ldots,k_j,\ldots,k_\nu) & = F(m,k_1,\ldots,k_j+1,\ldots,k_\nu),\\
	(QF) \left(m,k_1,\ldots,k_\nu \right) &= q^n F \left(m,k_1,\ldots,k_\nu \right),\\
	(Q_{(0)}F)(m,k_1,\ldots,k_\nu) &= q^m F(m,k_1,\ldots,k_\nu),\\
	(Q_jF)(m,k_1,\ldots,k_j,\ldots,k_\nu) & = q^{k_j}F(m,k_1,\ldots,k_j,\ldots,k_\nu),
\end{align*}
and let $\varepsilon$ be an evaluation map at $q=1$. Note that $E_{(0)}=E^2$ since $n=2m+1$.
Then, we verify the following proposition:
\begin{prop} \label{prop:intro1}
	The system of equations
	\begin{align*}
		\exp \left(w_j \frac{\partial \Phi _F}{\partial w_j} \right)&= 1, \quad(j = 1,\ldots,\nu)\\
		\exp \left(\alpha \frac{\partial \Phi _F}{\partial \alpha}\right) &= l^2
	\end{align*}
	coincides with
	\begin{align*}
		\varepsilon \left.\frac{E_j F}{F} \right|_{\substack{q^{k_{i}}=Q_{i} \\ q^m=Q_{(0)}}}&= 1, \quad(j = 1,\ldots,\nu)\\
		\varepsilon \left.\frac{E_{(0)} F}{F}\right|_{\substack{q^{k_{i}}=Q_{i} \\ q^m=Q_{(0)}}} &= E^2,
	\end{align*}
under the correspondences $w_j = Q_j$, $\alpha = Q_{(0)}$, and $ l = E$.
\end{prop}
As a corollary, we can obtain the polynomial $a_F(l,\alpha)$ from $q$-diffrences of the summand.
We can obtain an annihilating polynomial of $J_K(n)$ from the summand $F(m,k_1,\ldots,k_\nu)$ by creative telescoping (see Section \ref{sec:AJ}).
Through this method, we obtain the polynomial $P_F(E,Q,E_1,\ldots,E_\nu)$ in
\[
	\langle SE-R,\ S_1E_1-R_1,\ldots,\ S_\nu E_\nu-R_\nu \rangle \cap \mathbb{Q}[q,Q]\langle E,E_1,\ldots,E_\nu \rangle,
\]
where $R,\ S,\ R_j$, and $S_j$, with $j=1,\ldots,\nu$, are polynomials in $\mathbb{Q}[q,Q,Q_1,\ldots,Q_\nu]$ given by
\[
	\left. \frac{E_{j}F}{F}\right|_{\substack{q^{k_{i}}=Q_{i} \\ q^m=Q_{(0)}}} = \frac{R_{j}}{S_{j}},\text{ and }
	\left. \frac{EF}{F}\right|_{\substack{q^{k_{i}}=Q_{i} \\ q^m=Q_{(0)}}} = \frac{R}{S}.
\]
We put $P_F^0(E,Q)=P_F(E,Q,1,\ldots,1)$.
Comparing this process with the definition of $a_F(l,\alpha)$, we have the following theorem:
\begin{thm} \label{thm:intro2}
	The polynomial $a_F(l,\alpha)$ is a factor of $\varepsilon P_F^0(l,\alpha^2)$ for any solvable summand $F(m,k_1,\ldots,k_\nu)$.
\end{thm}
This paper is organized as follows:
In Section \ref{sec:pfApoly}, we review the colored Jones polynomial, its potential function and the $A$-polynomial.
Moreover we introduce the polynomial $a_D(l,\alpha)$ and see its property.
In Section \ref{sec:AJ}, we briefly look over the $A_q$-polynomial and creative telescoping.
In Section \ref{sec:polypf}, we introduce the polynomial $a_F(l,\alpha)$.
In Section \ref{sec:comp}, we compare those themes and verify Proposition \ref{prop:intro1} and Theorem \ref{thm:intro2}.
In the Appendix, we give the example of calculations for the figure-eight knot.
\par
\noindent \textit{Acknowledgments.} The author is grateful to Jun Murakami and Seokbeom Yoon for their helpful comments.

\section{Potential Functions and $A$-polynomials} \label{sec:pfApoly}
\subsection{Colored Jones polynomials}
\begin{figure}[htb]
	\centering
	\includegraphics[scale=0.8]{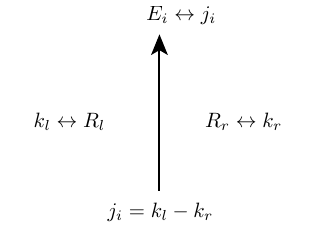}
	\caption{Let $j_i$ be an index assigned to an edge $E_i$,
	$k_l$ be an index assigned to the left region $R_l$ of the edge $E_i$,
	and$k_r$ be an index assigned to the right region $R_r$ of the edge $E_i$.
	These indices satisfy $j_i = k_l -k_r$.}
	\label{fig:regionedge}
\end{figure}

\begin{figure}[htb]
	\centering
	\includegraphics[scale=0.8]{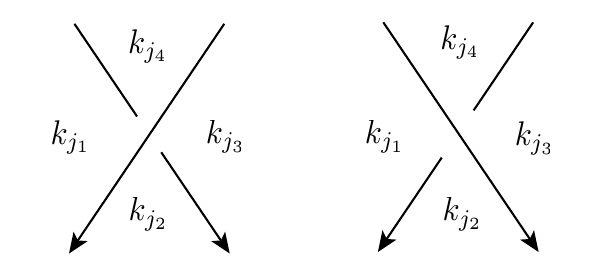}
	\caption{The indices labeled to the regions around a crossing.}
	\label{fig:regionind}
\end{figure}

In this section we review the definition of the colored Jones polynomial using the $R$-matrix.
First of all, we review the enhanced Yang-Baxter operator \cite{Li}.
\begin{defi} \label{def:eYB}
Let $V$ be a finitely generated free module over a commutative ring $\mathbb{K}$,
let $R:V \otimes V \to V \otimes V$ and $\mu:V \to V$ be automorphisms, $\id_V$ be the identity map,
and let $\alpha$ and $\beta$ be fixed units in $\mathbb{K}$.
We call $(R,\mu,\alpha,\beta)$ an enhanced Yang-Baxter operator if it satisfies the following conditions:
\begin{enumerate}
	\item $(R \otimes \id_V) \circ (\id_V \otimes R) \circ (R \otimes \id_V) = (\id_V \otimes R) \circ (R \otimes \id_V) \circ (\id_V \otimes R)$,
	\item $R \circ (\mu \otimes \mu) = (\mu \otimes \mu) \circ R$,
	\item $\tr _2(R^{\pm 1} \circ (\id_V \otimes \mu)) = \alpha^{\pm 1}\beta \id_V$.
\end{enumerate}
Here, $\tr_n: \mathrm{End}(V^{\otimes n}) \to \mathrm{End}(V^{\otimes(n-1)})$ is a map defined as follows:
For $f \in \mathrm{End}(V^{\otimes n})$, we put
\[
	f(e_{i_1} \otimes \cdots \otimes e_{i_n}) = \sum _{1 \leq j_1,\ldots,j_n \leq n}f^{j_1 \cdots j_n}_{i_1 \cdots i_n} e_{j_1} \otimes \cdots \otimes e_{j_n},
\]
where $\{e_1,\ldots,e_n\}$ is a besis of $V$.
Then, $\tr_n(f)$ is
\[
	\tr _n(f)(e_{i_1} \otimes \cdots \otimes e_{i_{n-1}}) = \sum _{1 \leq j_1,\ldots,j_{n-1},i_n \leq n}f^{j_1 \cdots j_{n-1}i_n}_{i_1 \cdots i_n} e_{j_1} \otimes \cdots \otimes e_{j_{n-1}}.
\]
\end{defi}
We can construct a link invariant from the enhanced Yang-Baxter operator.
\begin{prop}
Let $L$ be a link, let $b$ be a $k$-braid whose closure is $L$,
and let $(R,\mu,\alpha,\beta)$ be an enhanced Yang-Baxter operator.
Then,
\[
	T(L)=\alpha^{-w(b)}\beta^{-k} \tr(\phi(b) \circ \mu^{\otimes k})
\]
is a well-defined link invariant.
Here, $w$ is a group homomorphism from the braid group $B_k$ to the additive group $\mathbb{Z}$
called a writhe defined by $w(\sigma_i)=1$.
\end{prop}
We can construct enhanced Yang-Baxter operators from quantum enveloping algebras and their representations.
The $R$-matrix for the colored Jones polynomial (\cite{KM}) derives from the algebra $A_r,\ (r \in \mathbb{Z}_{>1})$
generated by $X,\ Y,\ K,\ \overline{K}$ with the relations
\begin{align*}
	\overline{K}=K^{-1},\quad KX=sXK,\quad &KY=s^{-1}YK,\quad XY-YX=\frac{K^2-\overline{K}^2}{s-s^{-1}},\\
	&X^r=Y^r=0,\quad K^{4r}=1,
\end{align*}
where $s=e^{\frac{\pi \iu}{r}}$, and its $n$-dimensional representation given by
\begin{align*}
	Xe_i &= [m+i+1]_s e_{i+1}, \\
	Ye_i &= [m-i+1]_s e_{i-1}, \\
	Ke_i &= s^{i}e_i,
\end{align*}
where $m$ is a half-integer satisfying $n=2m+1$,
$\{e_{m},e_{m-1},\ldots,e_{-m}\}$ is a basis of an $n$-dimensional complex vector space $V$,
and $[k]_s$ for an integer $k$ is a quantum integer defined by
\[
	[k]_s=\frac{s^k-s^{-k}}{s-s^{-1}}.
\]
Then, the $R$-matrix $R$ and the operator $\mu$ are given by
\begin{align*}%sign
	R(e_i \otimes e_j) = \sum _{k=0}^{\min \{m-i,\ m+j \}} & \frac{\{m-j+k\}!\{m+i+k\}!}{\{k\}!\{m-j\}!\{m+i\}!} \\
	& \times q^{ij-\frac{k(i-j)}{2}-\frac{k(k+1)}{4}}e_{j-k} \otimes e_{i+k},\\
	\mu(e_i)&=q^i e_i,\quad(i,j=-m,\ldots,m),
\end{align*}
where $s=q^{\frac{1}{2}}$, $\{k\} = q^{\frac{k}{2}}-q^{-\frac{k}{2}}$ for an integer $k$,
and $\{k\}!$ is a quantum factorial defined by
\[
	\{k\}!=\{k\}\{k-1\}\cdots \{1\},\quad \{0\}! = 1.
\]
\begin{rem}
The inverse of the above $R$ is given by
\begin{align*} %\label{eq:Rmat}
	\begin{split}
		R^{-1}(e_i \otimes e_j)=\sum _{k=0}^{\min \{m+i,\ m-j \}} &(-1)^{k}\frac{\{m-i+k\}!\{m+j+k\}!}{\{k\}!\{m-i\}!\{m+j\}!} \\
		& \times q^{-ij-\frac{k(i-j)}{2}+\frac{k(k+1)}{4}}e_{j+k} \otimes e_{i-k},\\
	\end{split}
\end{align*}
\end{rem}
We can calculate the colored Jones polynomial for a knot $K$ as follows:
We first present the knot $K$ as the closure of a $p$-braid $b \in B_p$.
Let $D$ be the corresponding diagram.
Then, we obtain an operator from $V^{\otimes p}$ to $V^{\otimes p}$
by assigning the $R$-matrix $R$ to each positive crossing
and its inverse $R^{-1}$ to each negative crossing,
and inserting the operator $\mu$ at the end of each string.
Modifying the trace of the operator so that the result is invariant under the Reidemeister move I,
we obtain the colored Jones polynomial for $K$.
Therefore, the colored Jones polynomial is of the form
\begin{gather*}
	J_n(K;q) = \sum _{i_1,\ldots,i_p} F^D(m,i_1,\ldots,i_p),\\
	F^D(m,i_1,\ldots,i_p) = q^{-(m^2+m)w(D)} q^{M(i_1,\ldots,i_p)}\prod _{c:\text{crossings}}(R^{\sgn(c)})^{i_\alpha i_\beta}_{i_\gamma i_\delta},
\end{gather*}
where $i_1,\ldots,i_p$ are indices assigned to edges of the diagram $D$,
$M(i_1,\ldots,i_p)$ is a polynomial that derives from the operator $\mu$
and has at most degree $1$ with respect to each index,
and $(R^{\sgn(c)})^{i_\alpha i_\beta}_{i_\gamma i_\delta}$ is a coefficient of the $R$-matrix assigned to the crossing $c$ with a sign $\sgn(c)$.
Namely, we put
\[
	R^{\pm 1}(e_{i_\gamma} \otimes e_{i_\delta}) = \sum _{i_\alpha,i_\beta}(R^{\pm})^{i_\alpha i_\beta}_{i_\gamma i_\delta}e_{i_\alpha} \otimes e_{i_\beta}.
\]
The indices here are assigned to the edges of the diagram.
For convenience later, we change those indices to the ones assigned to the regions of the diagram
according to the rules shown in Figure \ref{fig:regionedge}.
This variable transformation yields $F^D(m,k_1,\ldots,k_\rho)$ from $F^D(m,i_1,\ldots,i_p)$.
%We can rewrite the coefficient of the $R$-matrix to $R^{\pm}(m,k_{j_1},k_{j_2},k_{j_3},k_{j_4})$,
Explicitly, the coefficients of the $R$-matrix
\[
	R^{\pm}(m,k_{j_1},k_{j_2},k_{j_3},k_{j_4}) = (R^{\pm})^{k_{j_4}-k_{j_1},k_{j_3}-k_{j_4}}_{k_{j_2}-k_{j_1},k_{j_3}-k_{j_2}},
\]
where $k_{j_1},\ldots,k_{j_4}$ are indices assigned to the regions around the crossing $c$ as shown in the Figure \ref{fig:regionind},
take the following form:
\begin{align} \label{eq:Fcp}
\begin{split}
	R^+&(m,k_{j_1},k_{j_2},k_{j_3},k_{j_4})\\
	&=(-1)^{\kappa}q^{-m \kappa+(k_{j_3}-k_{j_4})(k_{j_4}-k_{j_1})+ \frac{\kappa^2-\kappa}{2}}\frac{(q)_{m +k_{j_1}-k_{j_4}}(q)_{m +k_{j_3}-k_{j_4}}}{(q)_{\kappa}(q)_{m+k_{j_2}-k_{j_3}}(q)_{m+k_{j_2}-k_{j_1}}},
\end{split}
\end{align}
\begin{align} \label{eq:Fcm}
\begin{split}
	R^- &(m,k_{j_1},k_{j_2},k_{j_3},k_{j_4})\\
	&=q^{m \kappa-(k_{j_2}-k_{j_1})(k_{j_3}-k_{j_2})}\frac{(q)_{m +k_{j_4}-k_{j_3}}(q)_{m +k_{j_4}-k_{j_1}}}{(q)_{-\kappa}(q)_{m+k_{j_1}-k_{j_2}}(q)_{m+k_{j_3}-k_{j_2}}},\\
\end{split}
\end{align}
where $ \kappa = k_{j_1}+k_{j_3}-k_{j_2}-k_{j_4}$, and $(q)_k =(1-q)(1-q^2) \cdots (1-q^k)$ for an integer $k$.
Here, we use the equality
\[
	\{k\}! = (-1)^k q^{-\frac{k(k+1)}{4}}(q)_k.
\]
The summand $F^D(m,k_1,\ldots,k_\rho)$ of the colored Jones polynomial is obtained directly from the diagram $D$
by multiplying the coefficients of the $R$-matrix associated to all crossings of $D$ and the factor derived from the operator $\mu$.
The actual formula of the colored Jones polynomial is often simplified using,
for example, a $q$-analogue of the binomial theorem.
Generally, the colored Jones polynomial is of the form
\begin{align} \label{eq:summFofcjp}
\begin{split}
	J_n(K;q) &=  \sum _{k_1,\ldots,i_\nu} F(m,k_1,\ldots,k_\nu),\\
	F(m,k_1,\ldots,k_\nu)&=(-1)^{(\boldsymbol{e},\tilde{\bk})}q^{\frac{1}{2}\tilde{\bk}^\top A \tilde{\bk}+(\bb,\tilde{\bk})+c}\frac{\prod_i (q)_{(\bp_i,\tilde{\bk})+r_i}}{\prod_j (q)_{(\bp'_j,\tilde{\bk})+r'_j}},
\end{split}
\end{align}
where $\boldsymbol{e}$ is a $(\nu+1)$-vector in $\{0,1\}^{\nu+1} \subset \mathbb{Q}^{\nu+1}$,
$ \tilde{\bk}=(m,k_1,\ldots,k_\nu)$,
$A$ is a $(\nu+1)$-th symmetric matrix with integer entries,
$\bb$ is a vector in $\mathbb{Q}^{\nu+1}$, $c$ is an integer, $\bp_i$ and $\bp'_j$ are $(\nu+1)$-vectors in $\mathbb{Z}^{\nu+1}\subset \mathbb{Q}^{\nu+1}$,
$r_i$ and $r'_j$ are integers.

\subsection{Potential functions} \label{sec:pf}
We want to observe the limit of the colored Jones polynomial for a knot $K$ with a diagram $D$ evaluated at the root of unity.
Let us briefly recall some facts on the potential function in the knot case. See \cite{Sa, Yon} for details.
\begin{defi}
	Suppose that the asymptotic behavior of a certain quantity $Q_N$ for a sufficiently large $N$ is %can be written as
	\[
		Q_N \sim \int \cdots \int _\Omega P_N e^{\frac{N}{2 \pi \iu} \Phi(z_1,\ldots,z_\nu)}dz_1 \cdots dz_\nu,
	\]
	where $P_N$ grows at most polynomially and $ \Omega $ is a region in $\mathbb{C}^\nu$. We call this function $\Phi(z_1,\ldots,z_\nu)$ a potential function of $Q_N$.
\end{defi}
A potential function $\Phi_F(\alpha,w_1,\ldots,w_\nu)$ of the colored Jones polynomial
evaluated at the root of unity is obtained by approximating the summand
$F(m,k_1,\ldots,k_\nu)$ of the colored Jones polynomial in the equation \eqref{eq:summFofcjp}
with continuous functions.
For some integers $n_1,\ n_2$ and $N$, we put $x_1= \xi_N^{n_1}$ and $x_2= \xi_N^{n_2}$, where $\xi_N = e^{\frac{2 \pi \iu}{N}}$.
Then the equalities
\[
	e^{\frac{2 \pi \iu}{N}n_1n_2} = e^{\frac{N}{2 \pi \iu}\left(\frac{2 \pi \iu}{N}n_1\right)\left(\frac{2 \pi \iu}{N}n_2\right)}, \quad \text{and}\quad
	(-1)^{n_1} = e^{\frac{N}{2 \pi \iu}\left(\pi \iu \frac{2 \pi \iu}{N}n_1\right)}
\]
respectively implies the correspondences
\begin{gather*}
	\xi_N^{mk_i} \leftrightarrow \log \alpha \log w_i ,\quad \text{and}\quad (-1)^{m} \leftrightarrow \pi \iu \log \alpha,\\
	\xi_N^{k_ik_j} \leftrightarrow \log w_i \log w_j ,\quad \text{and}\quad (-1)^{k_i} \leftrightarrow \pi \iu \log w_i,.
\end{gather*}
Moreover, the direct calculation shows
\begin{align*}
		\log(\xi_N)_{k} &= \sum^{k}_{j=1}\log(1-e^{\frac{2 \pi j \iu }{N}}) \\
		&= N \left( \int^{\frac{{k}}{N}}_0 \log(1-e^{2 \pi \iu \theta})d \theta+o(1) \right)\\
		&=\frac{N}{2 \pi \iu}\left( \int^{\xi^{k}_N}_1 \frac{\log(1-x)}{x}dx +o(1) \right)\\
		&= \frac{N}{2 \pi \iu}\left( -\li(\xi^{k}_N)+\frac{\pi^2}{6}+o(1)\right)
\end{align*}
for an integer $k$, which implies the correspondence
\[
	(\xi_N)_{(\bp,\tilde{\bk})+r} \leftrightarrow -\li(\tilde{\bw}^{\bp})+\frac{\pi^2}{6},
\]
where $\tilde{\bw}^{\bp} = \alpha^{p_0}w_1^{p_1} \cdots w_\nu^{p_\nu}$, with $\bp=(p_0,p_1,\ldots,p_\nu)$,
$\li(z)$ is a dilogarithm function defined by
\[
	\li(z) = -\int^z_0 \frac{\log(1-x)}{x}dx.
\]
A potential function $\Phi _F(\alpha,w_1,\ldots,w_\nu)$
of the colored Jones polynomial for a knot $K$ evaluated at $\xi_N$ is the summation of all these terms.
\begin{rem}
	$\xi_N$ to the power of a polynomial with at most degree $1$ does not contribute to the potential function.
\end{rem}
We will review some facts on the potential function \cite{Sa}.
Let $K$ be a hyperbolic knot with a diagram $D$,
and let $F^D(m,k_1,\ldots,k_\rho)$ be the summand of the colored Jones polynomial
obtained directly from the formulation using the $R$-matrix
(namely, without using a formula such as a $q$-analogue of the binomial theorem).
We can construct the potential function from the summand $F^D(m,k_1,\ldots,k_\rho)$.
We put it $\Phi_D(\alpha,w_1,\ldots,w_\rho)$ instead of $\Phi_{F^D}(\alpha,w_1,\ldots,w_\rho)$.
This potential function essentially coincides with the generalized potential function
in \cite{Yon}\footnote{In \cite{Yon}, Yoon defined the generalized potential function and established the relationship between the gluing equation.
The author appreciates Yoon's valuable comment at the 18th East Asian Conference on Geometric Topology.}. %%Yoon
From the saddle point of the potential function $\Phi_D(\alpha,w_1,\ldots,w_\rho)$,
we can obtain a hyperbolic structure of the complement $M_K$ of the knot $K$ that is not necessarily complete.
This is because
\[
	\exp \left(w_{i}\frac{\partial \Phi_D}{\partial w_{i}}\right) = 1,\quad i=1,\ldots,\rho
\]
coincides with the gluing equation of an ideal triangulation of $M_K$.
Since we assume that $K$ is a hyperbolic knot,
the equations above have a solution that gives the complete hyperbolic structure.
Then, we choose the saddle point $(\sigma_1(\alpha),\ldots,\sigma_\rho(\alpha))$
so that the imaginary part $\im \Phi$ of the potential function $\Phi$ satisfies %(\boldsymbol{\alpha},w_1,\ldots,w_\nu)
\[
	\im \Phi_D(\alpha,\sigma_1(\alpha),\ldots,\sigma_\rho(\alpha))|_{\alpha=-1} = \Vol(K),
\]
where $\Vol(K)$ is the hyperbolic volume of the knot $K$.
Namely, $(\sigma_1(-1),\ldots,\sigma_\rho(-1))$ corresponds to the complete hyperbolic structure with the volume $\Vol(K)$,
and $(\sigma_1(\alpha),\ldots,\sigma_\rho(\alpha))$ corresponds to the one-parameter deformation of the complete hyperbolic structure.
Here, $\alpha^2$ equals the dilation component of the action of the meridian,
and the hyperbolic structure obtained from the saddle point $(\sigma_1(\alpha),\ldots,\sigma_\rho(\alpha))$ is the cone-manifold.
Varying the value of $\alpha$ alters the cone angle along the knot.
Moreover, the following equality holds:
\[
	\exp \left(\alpha \frac{\partial \Phi_D}{\partial \alpha}(\alpha,\sigma_1(\alpha),\ldots,\sigma_\rho(\alpha))\right)=l(\alpha)^2
\]
where $l(\alpha)^2$ is the dilation component of the action of the preferred longitude of the knot $K$.
In other words, if% the system of equations
\begin{equation} \label{eq:pfsaddle}
	\begin{dcases}
		\exp \left(w_i \frac{\partial \Phi_D}{\partial w_i} \right)&= 1, (i = 1,\ldots,\rho)\\
		\exp \left(\alpha \frac{\partial \Phi_D}{\partial \alpha}\right) &= l^2,
	\end{dcases}
\end{equation}
has a solution, then $M_K$ admits a corresponding hyperbolic structure.
Note that the equations 
\begin{equation} \label{eq:pfsaddlegen}
	\begin{dcases}
		\exp \left(w_i \frac{\partial \Phi_F}{\partial w_i} \right)&= 1, (i = 1,\ldots,\nu)\\
		\exp \left(\alpha \frac{\partial \Phi_F}{\partial \alpha}\right) &= l^2,
	\end{dcases}
\end{equation}
are equivalent to algebraic equations for any summand $F(m,k_1,\ldots,k_\nu)$ of
the colored Jones polynomial.
%%%%cite pf aj and rt

\subsection{$A$-polynomial}%cite
In this subsection, we review the definition of the $A$-polynomial following \cite{CL}.
Let $K$ be a knot, and let $M_K$ be its complement.
The set of all $SL(2;\mathbb{C})$-representations $R(M_K)$ is an affine algebraic variety.
Then, we restrict our attention to
\[
	R_U=\{\varrho \in R(M_K) \mid \varrho(\mu) \text{ and } \varrho(\lambda) \text{ are upper triangular}. \},
\]
where $\mu$ is the meridian, and $\lambda$ is the preferred longitude of $K$.
We can define the eigenvalue map $\xi:R_U \to \mathbb{C}^2 $ by $\xi(\varrho)=(l,\alpha)$, where
\[
	\varrho(\lambda)=\left(
	\begin{array}{cc}
		l&	\ast \\
		0&	l^{-1}
	\end{array}
	\right),\quad \varrho(\mu)=\left(
	\begin{array}{cc}
		\alpha&	\ast \\
		0&	\alpha^{-1}
	\end{array}
	\right).
\]
For an algebraic component $C$ of $R_U$, the Zariski closure $\overline{\xi(C)}$ of $\xi(C)$ is an algebraic subset of $\mathbb{C}^2$.
If $\overline{\xi(C)}$ is a curve, there exists a defining polynomial.
The $A$-polynomial of a knot $K$ is the product of all such defining polynomials.
Since the $A$-polynomial contains a factor $l-1$, we call
\[
	A'_K(l,\alpha) = \frac{A_K(l,\alpha)}{l-1}.
\]
the nonabelian $A$-polynomial \cite{LZ}.
In the previous subsection, we recalled that the saddle point equation
coincides with a necessary condition where the knot complement admits a hyperbolic structure.
In fact, the factor $A'_K(l,\alpha)$ for some knots can be obtained from the saddle point equation of the potential function
of the colored Jones polynomial by eliminating parameters except $l$ and $\alpha$ \cite{Hi,Y2}.

\section{$A_q$-polynomial and AJ conjecture} \label{sec:AJ}
In this section, we recall some facts on the $A_q$-polynomial following \cite{Ga,GL}. %cite cjp is q-hol, Takata
\subsection{$A_q$-polynomial}
For a knot $K$, its colored Jones polynomial $J_K(n)=J_n(K;q)$ has a nontrivial linear recurrence relation \cite{GL} %%cite
\begin{equation} \label{eq:cjprec}
	\sum^d_{j=0}c_j(q,q^n)J_K(n+j)=0,
\end{equation}
where $c_j$ is a polynomial with integer coefficients.
For a discrete function $f: \mathbb{N} \to \mathbb{Q}(q)$, we define operators $Q$ and $E$ by
\[
	(Qf)(n)=q^nf(n),\quad (Ef)(n)=f(n+1).
\]
Then, we can restate the recurrence relation \eqref{eq:cjprec} as
\[
	\left(\sum^d_{j=0}c_j(q,Q)E^j\right)J_K(n) = 0.
\]
This polynomial $\sum c_j(q,Q)E^j$ is an element in the non-commutative algebra
\[
	\A = \mathbb{Z}[q^{\pm 1}]\langle Q,E \rangle/(EQ = qQE).
\]
To define the $A_q$-polynomial, we have to localize the algebra $\A$.
Let $\sigma$ be the automorphism of the field $\mathbb{Q}(q,Q)$ given by
\[
	\sigma(f)(q,Q)=f(q,qQ).
\]
The Ore algebra $ \A_{\mathrm{loc}}=\mathbb{Q}(q,Q)[E,\sigma]$ is defined by
\[
	\A_{\mathrm{loc}}=\left\{\sum^{\infty}_{k=0}a_k E^k \mid a_k \in \mathbb{Q}(q,Q),\ a_k=0 \text{ for sufficiently large }k\right\},
\]
where the multiplication of monomials is given by $aE^i \cdot bE^j = a\sigma^i(b)E^{i+j}$.
The $A_q$-polynomial $A_q(K)(E,Q)$ is a generator of the recursion ideal of $J_K(n)$%more explanation?
\[
	I=\{P \in \A_{\mathrm{loc}} \mid PJ_K(n)=0\},
\]
with the following properties:
\begin{itemize}
	\item $A_q(K)$ has the smallest $E$-degree and lies in $\A$.
	\item $A_q(K)$ is of the form $A_q(K) = \sum _k a_k E^k$, where $a_k \in \mathbb{Z}[q,Q]$ are coprime.
\end{itemize}
Garoufalidis proposed the following conjecture on the $A_q$-polynomial:
\begin{conj}[the AJ conjecture \cite{Ga}] \label{conj:AJ}
	For any knot $K$, $A_K(l,\alpha)$ is equal to $ \varepsilon A_q(K)(l,\alpha^2)$
	up to multiplication by an element in $ \mathbb{Q}(\alpha)$,
	where $ \varepsilon $ is an evaluation map at $q=1$.
\end{conj}
\subsection{Computation of $A_q$-polynomial}
\begin{defi}
	For a discrete function
	\[
		F:\mathbb{Z}^{\nu+1} \ni (n,k_1,\ldots,k_\nu) \to F(n,k_1,\ldots,k_\nu) \in \mathbb{Z}[q^\pm],
	\]
	we define operators $Q$, $E$, $Q_i$ and $E_i$ by
	\begin{align*}
		(QF)(n,k_1,\ldots,k_\nu) &= q^n F(n,k_1,\ldots,k_\nu), \\
		(EF)(n,k_1,\ldots,k_\nu) &= F(n+1,k_1,\ldots,k_\nu), \\
		(Q_i F)(n,k_1,\ldots,k_\nu) &= q^{k_i}F(n,k_1,\ldots,k_\nu), \\
		(E_i F)(n,k_1,\ldots,k_\nu) &= F(n,k_1,\ldots,k_i+1,\ldots,k_\nu).
	\end{align*}
\end{defi}
These operators generate the algebra %$E_i$ and $Q_i$, with $i=1,\ldots,r$,
\[
	\mathbb{Q}[q,Q,Q_1,\ldots,Q_\nu]\langle E,E_1,\ldots,E_\nu \rangle
\]
with relations
\begin{align*}
	& Q_iQ_j = Q_j Q_i,\ E_iE_j = E_jE_i,\ E_iQ_i = qQ_iE_i,\\
	& E_iQ_j = Q_jE_i \text{ for } i \neq j \in \{0,\ldots,\nu\},
\end{align*}
where $Q_0 = Q$ and $E_0=E$.
Hereafter, we put $ \boldsymbol{k}=(k_1,\ldots,k_\nu)$.
\begin{defi}
	A discrete function $F(n,\boldsymbol{k})$ is called $q$-hypergeometric if
	$E_iF/F \in \mathbb{Q}(q,q^n,q^{n_1},\ldots,q^{n_\nu})$ holds for all $i \in \{0,\ldots,\nu\}$.
\end{defi}
We especially deal with a proper $q$-hypergeometric function.
\begin{defi}
	A $q$-hypergeometric discrete function $F(n,\boldsymbol{k})$ is called proper if it is of the form
	\[
		F(n,\boldsymbol{k})=\frac{\prod _s(A_s;q)_{a_s n + (\boldsymbol{b}_s, \boldsymbol{k}) + c_s}}{\prod _t(B_t;q)_{u_t n + (\boldsymbol{v}_t, \boldsymbol{k}) + w_t}}q^{A(n,\boldsymbol{k})}\xi^{\boldsymbol{k}},
	\]
	where $A_s,B_t \in \mathbb{Q}(q)$, $a_s,u_t$ are integers,
	$\boldsymbol{b}_s,\boldsymbol{v}_t$ are vectors of $\nu$ integers,
	$c_s,w_t$ are integers that may depend on parameters,
	$A(n,\boldsymbol{k})$ is a quadratic form,
	$\xi$ is an $r$ vector of elements in $\mathbb{Q}(q)$, and
	\[
		(A;q)_n = \prod^{n-1}_{i=0}(1-Aq^i).
	\]
\end{defi}
Proper $q$-hypergeometric functions satisfy the following theorem:
\begin{thm}[\cite{WZ}]% cite/ meaning of proper
	Every proper $q$-hypergeometric function $F(n,\boldsymbol{k})$ has a $\boldsymbol{k}$-free recurrence
	\[
		\sum_{(i,\boldsymbol{j}) \in S} \sigma _{i,\boldsymbol{j}}(q^n)F(n+i,\boldsymbol{k}+\boldsymbol{j})=0,
	\]
	where $\boldsymbol{j}=(j_1,\ldots,j_\nu)$ is a $\nu$-tuple of integers, and
	$S$ is a finite set.
\end{thm}
We put
\begin{equation} \label{eq:annpolyF}
	P=P(E,Q,E_1,\ldots,E_\nu) = \sum_{(i,\boldsymbol{j}) \in S} \sigma _{i,\boldsymbol{j}}(Q)E^i E^{\boldsymbol{j}},
\end{equation}
where $E^{\boldsymbol{j}} = E_{1}^{j_1} \cdots E_{\nu}^{j_\nu}$.
$E$ and $Q$ are operators satisfying $(EF)(n,\boldsymbol{k})=F(n+1,\boldsymbol{k}),\ (QF)(n,\boldsymbol{k})=q^nF(n,\boldsymbol{k})$, and
this polynomial $P$ is an element in $\mathrm{Ann}(F) \cap \mathbb{Q}[q,Q]\langle E,E_k\rangle $.
Therefore, we would be able to obtain such $P$ by eliminating $Q_1,\ldots,Q_r$ from the system of equations
\[
	\begin{dcases}
		S_i E_i-R_i = 0 \quad (i=1,\ldots,r),\\
		SE-R=0,
	\end{dcases}
\]
where $EF/F = (R/S) \mid _{Q=q^n,\ Q_j = q^{n_j}}$.
The expansion of $P$ around $E_i=1$, with $i=1,\ldots,\nu$, is
\[
	P_0(E,Q)+\sum^\nu_{i=1}(E_i-1)R_i(E,Q,E_1,\ldots,E_\nu),
\]
where $P_0(E,Q)=P(E,Q,1,\ldots,1)$, and $R_i$ is a polynomial in $\mathbb{Q}[q,Q]\langle E,E_{\boldsymbol{k}} \rangle $.
Here, $E_{\boldsymbol{k}} = (E_1,\ldots,E_\nu)$.
Putting $G_i=R_iF$, we have
\[
	P_0(E,Q)F(n,\boldsymbol{k}) + \sum^\nu_{i=1}(G_i(n,k_1,\ldots,k_i+1,\ldots,k_\nu)-G_i(n,k_1,\ldots,k_\nu))=0.
\]
Summing up this equality, %$P_0(E,Q)G(n)=\mathrm{error}(n)$,
we verify that $G(n):=\sum _{\boldsymbol{k}} F(n,\boldsymbol{k})$ satisfies
\begin{align*}
	P_0(E,Q)G(n) &= \sum _{k_2,\ldots,k_\nu} (G_1(n,K_1,k_2,\ldots,k_\nu)-G_1(n,k^{0}_1,k_2,\ldots,k_\nu))\\
	+ \cdots &+ \sum _{\boldsymbol{k} \text{ except }k_i} (G_i(n,k_1,\ldots,K_i,\ldots,k_\nu)-G_i(n,k_1,\ldots,k^{0}_i,\ldots,k_\nu))\\
	+ \cdots &+ \sum _{k_1,\ldots,k_{\nu-1}} (G_\nu(n,k_1,\ldots,k_{\nu-1},K_\nu)-G_\nu(n,k_1,\ldots,k_{\nu-1},k^{0}_\nu)),
\end{align*}
where $K_i$ and $k^{0}_i$ are fixed parameters.
Since $F$ is a proper $q$-hypergeometric function, $G_i$ is a sum of proper $q$-hypergeometric functions.
That means $P_0(E,Q)G(n)$ is a sum of multisums of proper $q$-hypergeometric functions with one variable less.
Repeating this process, we obtain a polynomial $P_1(E,Q)$ such that
\[
	P_1(E,Q)P_0(E,Q)G(n)=0.
\]
\subsection{$A_q$-polynomial and eliminations}
The annihilating polynomial \eqref{eq:annpolyF} of the summand $F(n,\boldsymbol{k})$ is an element in
$\mathrm{Ann}(F) \cap \mathbb{Q}[q,Q]\langle E,E_{\boldsymbol{k}}\rangle $.
Moreover, defining the polynomials $R,S,R_i$, and $S_i \in \mathbb{Z}[q,Q_1,\ldots,Q_\nu]$ by
\[
	\frac{EF}{F} = \left. \frac{R}{S} \right|_{\substack{Q=q^n,\\ Q_j = q^{k_j}}},\quad \frac{E_iF}{F} =\left. \frac{R_i}{S_i} \right|_{\substack{Q=q^n,\\ Q_j = q^{k_j}}},
\]
it is known that the annihilating ideal of $F$ is generated by
\[
	\{SE-R\} \cup \{S_iE_i-R_i | i=1,\ldots,\nu\} \subset \mathbb{Q}[q,Q,Q_{\boldsymbol{k}}]\langle E,E_{\boldsymbol{k}} \rangle.
\]
Namely, in the process of obtaining an annihilating poynomial $P(E,Q,E_1,\ldots,E_\nu)$ of $F(n,\boldsymbol{k})$,
$Q_1,\ldots,Q_\nu$ are eliminated from
\[
	\begin{dcases}
		S_i E_i-R_i = 0 \quad (i=1,\ldots,\nu),\\
		SE-R=0.
	\end{dcases}
\]
In the Appendix \ref{ssec:annpoly41}, we calculate an annihilating polynomial of the colored Jones polynomial
for the figure-eight knot and demonstrate how the elimination is performed.
In this calculation, we use the Ore condition for $E+qQ$ and $1+QE$.
Though it appears ad hoc, it is the simplest method of elimination.
From the observation above, $ \varepsilon P_0(E,Q)$ is in
\[
	\varepsilon P_0(E,Q) \in \langle \varepsilon (SE-R),\ \varepsilon (S_1-R_1),\ldots,\ \varepsilon (S_\nu-R_\nu)\rangle \cap \mathbb{Q}[Q,E].
\]

\section{Polynomial obtained from a potential function} \label{sec:polypf}
In this section, we define a polynomial obtained from the potential function $\Phi_F(\alpha,w_1,\ldots,w_\nu)$.
Recall that the saddle point equation \eqref{eq:pfsaddlegen} of the potential function is equivalent to algebraic equations.
We view those equations as $\mathbb{Q}(\alpha)$-coefficient algebraic equations.
Namely, the saddle point equation is equivalent to the equations of the form
\begin{equation} \label{eq:peq0}
	p_0(l,w_1,\ldots,w_\nu) = 0,\ p_1(w_1,\ldots,w_\nu) = 0,\ldots,\ p_{\nu}(w_1,\ldots,w_\nu) = 0,
\end{equation}
where $p_0(l,w_1,\ldots,w_\nu) \in \mathbb{Q}(\alpha)[l,w_1,\ldots,w_\nu]$, and
$p_1(w_1,\ldots,w_\nu),\ldots,p_\nu(w_1,\ldots,w_\nu)$ are polynomials in $\mathbb{Q}(\alpha)[w_1,\ldots,w_\nu]$.
We suppose that those equations have solutions.
\begin{defi}
	We say that a summand $F(m,k_1,\ldots,k_\nu)$ is solvable if the system of equations \eqref{eq:pfsaddlegen}
	(or equivalently, the system of equations \eqref{eq:peq0}) has solutions.
\end{defi}
At least, the summand obtained directly from a knot diagram is solvable
because there exists a solution that gives a hyperbolic structure.
Then, the intersection
\[
	I = \langle p_0,p_1,\ldots,p_\nu \rangle \cap \mathbb{Q}(\alpha)[l]
\]
is an ideal in $\mathbb{Q}(\alpha)[l]$, where
$ \langle p_0,p_1,\ldots,p_\nu \rangle $ is the ideal in $\mathbb{Q}(\alpha)[l,w_1,\ldots,w_\nu]$
generated by $p_0,p_1,\ldots,p_\nu$.
Since $\mathbb{Q}(\alpha)[l]$ is a principal ideal domain, the ideal $I$ is generated by a single polynomial.
\begin{defi}
	For a solvable summand $F(m,k_1,\ldots,k_\nu)$, we define the polynomial $a_F(l,\alpha)$ as the generator of the ideal $I$ whose coefficients are in $\mathbb{Z}[\alpha]$ and coprime.
	The polynomial $a_F(l,\alpha)$ is defined up to multiplication by $ \pm 1$.
\end{defi}
Gukov \cite{Gu} describes the relationship between the colored Jones polynomial and
the $A$-polynomial from a mathematical physics perspective.
In the case of the summand $F^D$ obtained directly from a diagram $D$,
the system of equations \eqref{eq:pfsaddle} is a necessary condition where the knot complement admits a hyperbolic structure \cite{Sa}.
Then, we expect that it is related to the $A$-polynomial.
In fact, the nonabelian $A$-polynomials for some knots are obtained from the potential function of the colored Jones polynomial
by eliminating redundant variables,
though the formulae are not the ones directly obtained from the $R$-matrix formulation \cite{Hi,Y2}.
Therefore, the following question is reasonable:
\begin{ques} \label{q:aD}
Does the polynomial $a_F(l,\alpha)$ coincide with the nonabelian $A$-polynomial (up to multiplication by $ \pm 1$)
for any solvable summand $F(m,k_1,\ldots,k_\nu)$ of the colored Jones polynomial?
\end{ques}
In the Appendix \ref{ssec:Apoly41}, we calculate the $A$-polynomial for the figure-eight knot
from the potential function, for example.

\section{Comparison with the saddle point equation} \label{sec:comp}
In this section, we compare the saddle point equation of the potential function of the
colored Jones polynomial for a hyperbolic knot $K$.
First of all, we recall that the colored Jones polynomial is the summation of
\[
	F(m,k_1,\ldots,k_\nu)=(-1)^{(\boldsymbol{e},\tilde{\bk})}q^{\frac{1}{2}\tilde{\bk}^\top A \tilde{\bk}+(\bb,\tilde{\bk})+c}\frac{\prod_i (q)_{(\bp_i,\tilde{\bk})+r_i}}{\prod_j (q)_{(\bp'_j,\tilde{\bk})+r'_j}},
\]
where $\boldsymbol{e}$ is a $(\nu+1)$-vector in $\{0,1\}^{\nu+1}$,
$A$ is a $(\nu+1)$-th symmetric matrix with integer entries,
$\bb$ is a vector in $\mathbb{Q}^{\nu+1}$, $c$ is an integer, $\bp_i$ and $\bp'_j$ are $(\nu+1)$-vectors in $\mathbb{Z}^{\nu+1}$, $r_i$ and $r'_j$ are integers.
We will compare the derivatives of the potential function and the $q$-differences of the summand.
Hereafter, we suppose that the summand $F(m,k_1,\ldots,k_\nu)$ is solvable.
Let $E,\ E_{(0)},\ E_j,\ Q,\ Q_{(0)}$, and $Q_j$ be operators defined by
\begin{align*}
	(EF)\left(\frac{n-1}{2},k_1,\ldots,k_\nu \right) &= F \left(\frac{n}{2},k_1,\ldots,k_\nu \right),\\
	(E_{(0)}F)(m,k_1,\ldots,k_\nu) &= F(m+1,k_1,\ldots,k_\nu),\\
	(E_jF)(m,k_1,\ldots,k_j,\ldots,k_\nu) & = F(m,k_1,\ldots,k_j+1,\ldots,k_\nu),\\
	(QF) \left(m,k_1,\ldots,k_\nu \right) &= q^n F \left(m,k_1,\ldots,k_\nu \right),\\
	(Q_{(0)}F)(m,k_1,\ldots,k_\nu) &= q^m F(m,k_1,\ldots,k_\nu),\\
	(Q_jF)(m,k_1,\ldots,k_j,\ldots,k_\nu) & = q^{k_j}F(m,k_1,\ldots,k_j,\ldots,k_\nu),
\end{align*}
where $n=2m+1$, and let $\varepsilon$ be an evaluation map at $q=1$.
Then, the following theorem holds:
\begin{prop} \label{prop:comp}
	The system of equations
	\begin{align}
		\exp \left(w_j \frac{\partial \Phi _F}{\partial w_j} \right)&= 1, \quad(j = 1,\ldots,\nu) \label{eq:diffpfw}\\
		\exp \left(\alpha \frac{\partial \Phi _F}{\partial \alpha}\right) &= l^2 \label{eq:diffpfa}
	\end{align}
	coincides with
	\begin{align}
		\varepsilon \left.\frac{E_j F}{F} \right|_{\substack{q^{k_{i}}=Q_{i} \\ q^m=Q_{(0)}}}&= 1, \quad(j = 1,\ldots,\nu) \label{eq:qdiffsumEj}\\
		\varepsilon \left.\frac{E_{(0)} F}{F}\right|_{\substack{q^{k_{i}}=Q_{i} \\ q^m=Q_{(0)}}} &= E^2, \label{eq:qdiffsumEm}
	\end{align}
under the correspondences $w_j = Q_j$, $\alpha = Q_{(0)}$, and $ l = E$.
\end{prop}
\begin{proof}
We put $A=(a_{ij})_{0 \leq i,j \leq \nu}$, $\bb=(b_0,\ldots,b_\nu)^\top$, and $\bp=(p_0,\ldots,p_\nu)^\top$.
Moreover, we put $k_0 = m$ and $w_0 = \alpha $ in the following argument.
We compare $\exp \left(w_j \frac{\partial \Phi}{\partial w_j} \right)$ and $\varepsilon \frac{E_j F}{F}$
with respect to the following three kinds of factors:
\[
	(-1)^{k_j},\quad q^{\frac{1}{2}\tilde{\bk}^\top A \tilde{\bk}+(\bb,\tilde{\bk})+c},\text{ and } (q)_{(\bp,\tilde{\bk})+r}^{\pm 1}
\]
For $(-1)^{k_j}$, the corresponding term of the potential function is $ \pi \iu \log w_j$.
Therefore, its derivative satisfies
\[
	\exp \left(w_j\frac{\partial}{\partial w_j}\pi \iu \log w_j\right) = -1.
\]
On the other hand,
\[
	\varepsilon \left.\frac{E_j (-1)^{k_j}}{(-1)^{k_j}} \right|_{\substack{q^{k_{i}}=Q_{i} \\ q^m=Q_{(0)}}}=-1.
\]
For $q^{\frac{1}{2}\tilde{\bk}^\top A \tilde{\bk}+(\bb,\tilde{\bk})+c}$, it is enough to consider the term
\[
	f(m,k_1,\ldots,k_\nu)=q^{\frac{1}{2}a_{jj}k_j^2+b_jk_j + \sum_{i \neq j}a_{ij}k_ik_j + c}.
\]
The corresponding term of the potential function is
\[
	\varphi(\alpha,w_1,\ldots,w_\nu)=\frac{a_{jj}}{2}(\log w_j)^2 + \sum_{i \neq j}a_{ij}\log w_i \log w_j.
\]
Therefore, its derivative satisfies
\[
	\exp \left(w_j\frac{\partial \varphi}{\partial w_j}\right) = \alpha^{a_{0j}}w_1^{a_{1j}} \cdots w_\nu^{a_{\nu j}}
\]
On the other hand,
\begin{align*}
	\varepsilon \left.\frac{E_j f}{f} \right|_{\substack{q^{k_{i}}=Q_{i} \\ q^m=Q_{(0)}}} &= \varepsilon (q^{\frac{a_{jj}}{2}+b_j + \sum a_{ij}k_i})|_{\substack{q^{k_{i}}=Q_{i} \\ q^m=Q_{(0)}}}\\
	&=Q_{(0)}^{a_{0j}}Q_1^{a_{1j}} \cdots Q_\nu^{a_{\nu j}}.
\end{align*}
Note that $\frac{a_{jj}}{2}+b_j \in \mathbb{Z}$ since $F_D(m,k_1,\ldots,k_\nu)$ is a $q$-hypergeometric function.
For $(q)_{(\bp,\tilde{\bk})+r}$, the corresponding term of the potential function is
\[
	\psi(\alpha,w_1,\ldots,w_\nu) = -\li(\tilde{\bw}^{\bp})+\frac{\pi^2}{6}.
\]
Therefore, its derivative is
\[
	\exp \left(w_j\frac{\partial \varphi}{\partial w_j}\right) = (1-\tilde{\bw}^{\bp})^{p_j}.
\]
On the other hand, if $p_j > 0$,
\begin{align*}
	\varepsilon \left.\frac{E_j (q)_{(\bp,\tilde{\bk})+r}}{(q)_{(\bp,\tilde{\bk})+r}} \right|_{\substack{q^{k_{i}}=Q_{i} \\ q^m=Q_{(0)}}} &=
	\varepsilon (1-q^{r+1}Q_{(0)}^{p_0}\cdots Q_\nu^{p_\nu}) \cdots (1-q^{r+p_j}Q_{(0)}^{p_0}\cdots Q_\nu^{p_\nu})\\
	&=(1-Q_{(0)}^{p_0}\cdots Q_\nu^{p_\nu})^{p_j}.
\end{align*}
If $p_j < 0$,
\begin{align*}
	\varepsilon \left.\frac{E_j (q)_{(\bp,\tilde{\bk})+r}}{(q)_{(\bp,\tilde{\bk})+r}} \right|_{\substack{q^{k_{i}}=Q_{i} \\ q^m=Q_{(0)}}} &=
	\varepsilon \frac{1}{(1-q^{p_k+r+1}Q_{(0)}^{p_0}\cdots Q_\nu^{p_\nu}) \cdots (1-q^{r}Q_{(0)}^{p_0}\cdots Q_\nu^{p_\nu})}\\
	&=(1-Q_{(0)}^{p_0}\cdots Q_\nu^{p_\nu})^{p_j}.
\end{align*}
We can apply a similar consideration to the factor of the form $(q)_{(\bp,\tilde{\bk})+r}^{-1}$.
Therefore, the systems of equations \eqref{eq:diffpfw} and the equation \eqref{eq:diffpfa} respectively coincide with
the systems of equations \eqref{eq:qdiffsumEj} and the equation \eqref{eq:qdiffsumEm}
under the correspondences $w_j = Q_j$, $\alpha = Q_{(0)}$, and $ l = E$.
\end{proof}
We define the polynomials $R_{j}$ and $S_{j}$, with $j=1,\ldots,\nu$, by
\[
	\left. \frac{E_{j}F}{F}\right|_{\substack{q^{k_{i}}=Q_{i} \\ q^m=Q_{(0)}}} = \frac{R_{j}}{S_{j}}.
\]
The system of equations $ \varepsilon(S_j E_j -R_j) = 0$, with $j=1,\ldots,\nu$, under the substitution $E_j = 1$,
which means
\[
	\varepsilon(S_j-R_j) = 0, \quad (j=1,\ldots,\nu),
\]
is equivalent to the system of equations \eqref{eq:qdiffsumEj}.
This corresponds to the system of equations \eqref{eq:diffpfw} by Proposition \ref{prop:comp}.
We also define polynomials $R$ and $S$ by
\[
	\left. \frac{EF}{F}\right|_{\substack{q^{k_{i}}=Q_{i} \\ q^m=Q_{(0)}}} = \frac{R}{S}.
\]
Note that $E_{(0)}=E^2$ holds because $n=2m+1$.
Then,
\[
	\left. \varepsilon \frac{E_{(0)}F}{F}\right|_{\substack{q^{k_{i}}=Q_{i} \\ q^m=Q_{(0)}}} = \left. \varepsilon \frac{E^2 F}{F}\right|_{\substack{q^{k_{i}}=Q_{i} \\ q^m=Q_{(0)}}}=\left. \varepsilon \frac{E^2 F}{EF}\frac{EF}{F}\right|_{\substack{q^{k_{i}}=Q_{i} \\ q^m=Q_{(0)}}} = \varepsilon \frac{R^2}{S^2}
\]
holds. The equation $ \varepsilon(S^2 E_{(0)} -R^2)= \varepsilon(S^2 E^2 -R^2) = 0$ is equivalent to
the equation \eqref{eq:qdiffsumEm}.
This corresponds to the equation \eqref{eq:diffpfa} by Proposition \ref{prop:comp}.
\begin{cor}
	We can obtain the polynomial $a_F(l,\alpha)$ by taking a generator of the ideal
	\[
		\langle \varepsilon(S_1-R_1),\ \ldots,\varepsilon(S_\nu-R_\nu),\ \varepsilon(S^2 E^2 -R^2)\rangle \cap \mathbb{Q}(Q_{(0)})[E]
	\]
	in $\mathbb{Q}(Q_{(0)})[E]$ whose coefficients are in $\mathbb{Z}[Q_{(0)}]$ and coprime, and substituting $E=l$ and $Q_{(0)}=\alpha $.
\end{cor}
Note that the system of equations
\[
	\begin{dcases}
		S_i E_i-R_i = 0 \quad (i=1,\ldots,\nu),\\
		S^2E^2-R^2=0.
	\end{dcases}
\]
yields an annihilating polynomial of the summand in $E,Q_{(0)},E_1,\ldots,E_\nu$,
but we can obtain the one $P_F(E,Q,E_1,\ldots,E_\nu)$ from it because $Q=qQ_{(0)}^2$.
Note also that $Q_{(0)}$ corresponds to the eigenvalue $\alpha$ of the meridian.
This explains the substitution of $m^2$ for the $A_q$-polynomial in the statement of the AJ conjecture.
When all indices $\boldsymbol{k}$ vanish after finite times of creative telescoping,
we obtain an inhomogeneous recurrence relation
\begin{equation} \label{eq:inhomJK}
	\tilde{P}_F(E,Q)J_K(n) + f(q,q^n) = 0,
\end{equation}
where $\tilde{P}_F(E,Q) \in \A_{\mathrm{loc}}$, and $f(q,q^n) \in \mathbb{Q}(q,q^n)$.
Since $f(q,q^n)$ can be canceled by left multiplication of $(E-1) \cdot f(q,Q)^{-1}$,
we obtain homogeneous recurrence relation
\[
	(E-1) \cdot \frac{1}{f(q,Q)} \cdot \tilde{P}_F(E,Q)J_K(n) = 0.
\]
This would support the AJ conjecture.
Note that $P_F^0(E,Q)= P_F(E,Q,1,\ldots,1)$ is a factor of $\tilde{P}_F(E,Q)$ by the procedure of creative telescoping.
Moreover, recall that
\begin{gather*}
	\varepsilon P_F^0(E,Q) \in \langle \varepsilon (S^2E^2-R^2),\varepsilon (S_1-R_1),\ldots,\varepsilon (S_\nu-R_\nu)\rangle \cap \mathbb{Q}[Q,E] \\
	\subset \langle \varepsilon (S^2E^2-R^2),\varepsilon (S_1-R_1),\ldots,\varepsilon (S_\nu-R_\nu)\rangle \cap \mathbb{Q}(Q)[E].
\end{gather*}
This implies the following:
\begin{thm}
	The polynomial $a_F(l,\alpha)$ is a factor of $\varepsilon P_F^0(l,\alpha^2)$
	for any solvable summand $F(m,k_1,\ldots,k_\nu)$.
\end{thm}
\begin{rem} \label{rem:approx}
We can view Proposition \ref{prop:comp} as follows:
The potential function is obtained from the summand of the colored Jones polynomial by approximating it with continuous functions.
Therefore, for a sufficiently large integer $N$,
\[
	F(m,\boldsymbol{k})|_{q=\xi_N} \sim P_N \exp \left(\frac{N}{2 \pi \iu}\Phi _F(\xi_N^m,\xi_N^{k_1},\ldots,\xi_N^{k_\nu}) \right).
\]%where $\xi_N=e^{\frac{2 \pi \iu}{N}}$. \substack{\alpha =\xi_N^m \\ w_i=\xi_N^{k_i}}
Then,
\begin{align*}
	\left. \frac{E_{(0)}F}{F} \right|_{q=\xi_N} &\sim \frac{\exp \left(\frac{N}{2 \pi \iu}\Phi _F(\xi_N \cdot \xi_N^m,\xi_N^{k_1},\ldots,\xi_N^{k_\nu})\right)}{\exp \left(\frac{N}{2 \pi \iu}\Phi _F(\xi_N^m,\xi_N^{k_1},\ldots,\xi_N^{k_\nu})\right)}\\
	&\sim \exp \left(\xi_N^m \frac{\partial \Phi _F}{\partial \alpha}(\xi_N^m,\xi_N^{k_1},\ldots,\xi_N^{k_\nu})\right).
\end{align*}
\end{rem}
%we would be able to obtain the factor of the $A$-polynomial $A_K(l,\alpha)$ for $K$ corresponding to nonabelian representations
%from $P_K^0(E,Q)$.

\appendix
\section{Example calculation for the figure-eight knot}
Let us observe the process above with the figure-eight knot.
The colored Jones polynomial for the figure-eight knot is \cite{Ha}
	\[
		J_n(4_1,q) = \sum^{n-1}_{i=0}F(n,i),
	\]
	where
	\[
		F(n,i) = \frac{1}{\{n\}} \frac{\{n+i\}!}{\{n-i-1\}!}.
	\]
Note that this formula is not the form itself obtained from the $R$-matrix.
The argument above, however, is still valid.
\subsection{Potential function and the $A$-polynomial} \label{ssec:Apoly41}
	The potential function $\Phi(\alpha,x)$ of $J_i(4_1,\xi_N)$ is (see \cite{Sa})
	\[
		\Phi(\alpha,x) = -2 \log \alpha \log x -\li(\alpha^2x)+\li(\alpha^2x^{-1}).
	\]
	The derivatives of $\Phi$ with $x$ and $\alpha$ are
	\begin{align*}
		x \frac{\partial \Phi}{\partial x} &= \log \alpha^{-2}(1-\alpha^2 x)(1-\alpha^2 x^{-1}),\\
		\alpha \frac{\partial \Phi}{\partial \alpha} &= 2\log(1-\alpha^2 x)(x-\alpha^2)^{-1}.
	\end{align*}
	Noting that if $\rho(\lambda)$ is of the form
	\[
		\rho(\lambda)=\left(
		\begin{array}{cc}
			l&	\ast \\
			0&	l^{-1}
		\end{array}
		\right)
	\]
	the action of $\lambda$ is $z \mapsto l^2 z + \ast $, especially its dilation component is equal to $l^2$, we put
	\[
		(1-\alpha^2 x)^2(x-\alpha^2)^{-2}=l^2
	\]
	From% the equations
	\[
		\begin{dcases}
			\alpha^{-2}(1-\alpha^2 x)(1-\alpha^2 x^{-1})=1,\\
			(1-\alpha^2 x)(x-\alpha^2)^{-1} =l,
		\end{dcases}
	\]%% needs some comments (PSL(2;C))
	we obtain the factor of the (nonabelian) A-polynomial of the figure-eight knot
	\begin{equation} \label{eq:A41}
		\alpha^4 l^2-l+\alpha^2 l + 2 \alpha^4 l + \alpha^6 l - \alpha^8 l + \alpha^4
	\end{equation}
	by eliminating $x$.
\subsection{Annihilating polynomials of $J_n(4_1,q)$}\label{ssec:annpoly41}
The annihilating polynomial of $J(n)=J_n(4_1;q)$ is \cite{Ga}
{\small
\begin{align*}
	&\frac{q^{4}Q(-1+q^{3}Q)}{(q+q^{3}Q)(q-q^{6}Q^{2})}E^3\\
	&+\frac{(-q+q^{3}Q)(q^4+q^5Q-2q^{6}Q-q^{7}Q^{2}+q^{8}Q^{2}-q^{9}Q^{2}-2q^{10}Q^{3}+q^{11}Q^3+q^{12}Q^{4})}{q^{4}Q(q^2+q^{3}Q)(-q+q^{6}Q^{2})}E^2\\
	&-\frac{(q^2-q^{3}Q)(q^8-2q^{9}Q+q^{10}Q-q^{9}Q^{2}+q^{10}Q^{2}-q^{11}Q^{2}+q^{10}Q^{3}-2q^{11}Q^{3}+q^{12}Q^{4})}{q^{5}Q(q+q^{3}Q)(q^5-q^{6}Q^{2})}E\\
	&+\frac{q^{5}Q(-q^3+q^{3}Q)}{(q^2+q^{3}Q)(-q^5+q^{6}Q^{2})}.
\end{align*}
}
We can factorize this polynomial as $(E-1)\alpha(q,E,Q)(Q-1)$,
where $\alpha(q,E,Q)$ is
\[
	\frac{1}{1+qQ} \left\{ \frac{qQ}{1-q^3Q^2}E^2 + \left( \frac{1}{1-q^3Q^2} + \frac{1}{1-qQ^2}+qQ-1-\frac{1}{qQ} \right)E + \frac{qQ}{1-qQ^2}\right\}.
\]
$EF/F$ and $E_1F/F$ are
	\begin{align} \label{eq:def41}
	\begin{split}
		\frac{EF}{F} &= \frac{F(n+1,i)}{F(n,i)} = \frac{(1-q^n)(1-q^{n+1+i})}{(1-q^{n+1})(q^i-q^n)} \\
		\frac{E_1 F}{F} &= \frac{F(n,i+1)}{F(n,i)} = q^{-n}(1-q^{n+i+1})(1-q^{n-i-1}).
	\end{split}
	\end{align}
Substituting $Q=q^n$ and $Q_1=q^i$ into \eqref{eq:def41},
we have %the difference equation of $F(n,i)$%%%%
\begin{align}
	(E+qQ)Q_1(Q-1) = (1+QE)(Q-1), \label{eq:41E}\\
	q^2 Q_1^2Q+q Q_1(-Q^2+Q E_1-1)+Q=0. \label{eq:41E1}
\end{align}
From \eqref{eq:41E}, we have
\begin{equation} \label{eq:41Q1}
	(1+QE)Q_1^{-1}(Q-1)=(E+qQ)(Q-1)
\end{equation}
Multiplying \eqref{eq:41E1} by $q^{-1}Q_1^{-1}Q^{-1}(Q-1)$ from the left, we obtain
\begin{equation} \label{eq:Q1Q1inv}
	q Q_1(Q-1)+Q^{-1}(-Q^2+Q E_1-1)(Q-1)+q^{-1} Q_1^{-1}(Q-1)=0.
\end{equation}
Then, we multiply \eqref{eq:Q1Q1inv} by
\[
	X(q,E,Q)=\frac{qQ}{1-q^3Q^2}E^2 + \left( \frac{1}{1-q^3Q^2}+\frac{1}{1-qQ^2}-1\right)E+\frac{qQ}{1-qQ^2}
\]
from the left. This polynomial is factorized in two ways.
\begin{align*}
	X(q,E,Q)&=\left(\frac{qQ}{1-q^3Q^2}E+\frac{1}{1-qQ^2} \right)(E+qQ)\\
	&= \left(\frac{1}{1-q^3Q^2}E+\frac{qQ}{1-qQ^2} \right)(1+QE).
\end{align*}
Then, using \eqref{eq:41E} and \eqref{eq:41Q1}, we obtain the annihilating polynomial $P(E,Q,E_1)$ of $F(n,i)$
\begin{align*}
	&P(E,Q,E_1)=\left\{ \frac{qQ}{1-q^3Q^2}E_1E^2 \right.\\
	 &+ \left. \left( \frac{1}{1-q^3Q^2}E_1 + \frac{1}{1-qQ^2}E_1+qQ-E_1-\frac{1}{qQ} \right)E + \frac{qQ}{1-qQ^2}E_1\right\}(Q-1).
\end{align*}
The expansion of $P(E,Q,E_1)$ at $E_1=1$ is
\[
	P(E,Q,E_1) = P_0(E,Q) + (E_1-1)R(E,Q),
\]
where
\begin{equation}\label{eq:41P0}
	P_0(E,Q) = P(E,Q,1)=(1+qQ)\alpha(q,E,Q)(Q-1),
\end{equation}
and 
\[
	R(E,Q) = \left\{\frac{qQ}{1-q^3Q^2}E^2 + \left( \frac{1}{1-q^3Q^2} + \frac{1}{1-qQ^2}-1 \right)E + \frac{qQ}{1-qQ^2} \right\}(Q-1).
\]
Therefore, $P_0(E,Q)F$ is of the form
\[
	P_0(E,Q)F = c_2(q,q^n)F(n+2,i) + c_1(q,q^n)F(n+1,i) + c_0(q,q^n)F(n,i),
\]
where $c_k(q,q^n) \in \mathbb{Q}(q,q^n)$, with $k=0,1,2$.
Summing up this equality with $i$ running from $0$ to $n+1 (=n+2-1)$, we have
\[
	P_0(E,Q)J(n) = c_2(q,q^n)J(n+2) + c_1(q,q^n)J(n+1) + c_0(q,q^n)J(n).
\]
Note that $F(n,i)= 0$ when $i \geq n$.
Putting $G(n,i)=R(E,Q)F$, on the other hand, we have
\[
	\sum^{n+1}_{i=0}(E_1-1)G(n,i) = G(n,n+2)-G(n,0) = q^{n+1}+1.
\]
Therefore, we have the second order inhomogeneous recurrence relation \cite{Ga}
\[
	P_0(E,Q)J(n)+q^{n+1}+1 = 0
\]
Since $q^{n+1}+1$ is annihilated by
\[
	P_1(E,Q) = (E-1) \cdot \frac{1}{1+qQ},
\]
we have the third order homogeneous recurrence relation $P_1(E,Q)P_0(E,Q)J(n)= 0$. %%subsection from here?
\subsection{Comparison of the derivatives and the $q$-differences}
	Substituting $q=1$ into \eqref{eq:41E} and \eqref{eq:41E1},
	we have% the following system of equations:
	\begin{align}%arrange%
		(E+Q)Q_1(Q-1) = (1+QE)(Q-1), \label{eq:e41E}\\
		Q_1^2Q+ Q_1(-Q^2+Q-1)+Q=0. \label{eq:e41E1}
	\end{align}
	Here, the factor $(Q-1)$ in \eqref{eq:e41E} is canceled and we have
	\begin{equation}
		(E+Q) Q_1 = 1+QE. \label{eq:e41Ecan}
	\end{equation}
	Eliminating $Q_1$ from \eqref{eq:e41E1} and \eqref{eq:e41Ecan}, we obtain
	\begin{equation} \label{eq:e41P0}
		Q^2E^2-E+QE-2Q^2E + Q^3E-Q^4E+Q^2,
	\end{equation}
	which is equal to the polynomial \eqref{eq:A41} under the substitutions $Q= \alpha^2$ and $E=l$.
	The polynomial \eqref{eq:e41P0} is also equal to the one \eqref{eq:41P0} with $q$ evaluated at $1$
	\[
		\varepsilon P_0(E,Q) = \frac{1}{Q(1-Q^2)}(Q^2E^2-E+QE-2Q^2E + Q^3E-Q^4E+Q^2),
	\]
	up to multiplication by an element in $ \mathbb{Q}(Q)$.

\end{document}